\documentclass[twocolumn]{autart}
\usepackage[latin9]{inputenc}
\usepackage{float}
\usepackage{amsmath}
\usepackage{amssymb}
\usepackage{graphicx}
\usepackage{epstopdf}
\usepackage{algorithm}
\usepackage{algpseudocode}
\makeatletter

\providecommand{\tabularnewline}{\\}
\floatstyle{ruled}
\newfloat{algorithm}{tbp}{loa}
\providecommand{\algorithmname}{Algorithm}
\floatname{algorithm}{\protect\algorithmname}

\usepackage{cite}
\usepackage[subnum]{cases}
\usepackage{color}

\usepackage{graphicx}

\usepackage{macros}

\makeatother

\begin{document}
\begin{frontmatter} 

\title{Construction of Parametric Barrier Functions \\
for Dynamical Systems using Interval Analysis\thanksref{footnoteinfo}}


\thanks[footnoteinfo]{This research was partially supported by Labex DigiCosme (project
ANR-11-LABEX-0045-DIGICOSME) operated by ANR as part of the program
``Investissement d'Avenir'' Idex Paris-Saclay (ANR-11-IDEX-0003-02)
and by the ANR INS Project CAFEIN (ANR-12-INSE-0007).}

\author[ENSTA]{A. Djaballah}\ead{adel.djaballah@ensta-paristech.fr},
\author[ENSTA]{A. Chapoutot}\ead{alexandre.chapoutot@ensta-paristech.fr},
\author[Michel]{M. Kieffer}\ead{michel.kieffer@lss.supelec.fr},
\author[CEA]{O. Bouissou}\ead{olivier.bouissou@cea.fr}

\address[ENSTA]{U2IS, ENSTA ParisTech, Universit\'e Paris-Saclay, 828 bd des Mar\'echaux, 91762 Palaiseau Cedex}

\address[Michel]{L2S, CNRS, Supélec, Univ. Paris-Sud 91192 Gif-sur-Yvette Cedex and
Institut Universitaire de France 75005 Paris, }

\address[CEA]{CEA Saclay Nano-INNOV Institut~CARNOT, 91191 Gif-sur-Yvette Cedex}


\begin{abstract}
Recently, barrier certificates have been introduced to prove the safety
of continuous or hybrid dynamical systems. A barrier certificate needs
to exhibit some barrier function, which partitions the state space
in two subsets: the safe subset in which the state can be proved to
remain and the complementary subset containing some unsafe region.
This approach does not require any reachability analysis, but needs
the computation of a valid barrier function, which is difficult when
considering general nonlinear systems and barriers. This paper presents
a new approach for the construction of barrier functions for nonlinear
dynamical systems. The proposed technique searches for the parameters
of a parametric barrier function using interval analysis. Complex
dynamics can be considered without needing any relaxation of the constraints
to be satisfied by the barrier function.

\end{abstract}
\end{frontmatter}

\section{Introduction}

\label{sec:introduction} Formal verification aims at proving that a certain behavior or property
is fulfilled by a system. Verifying,~\emph{e.g.}, the safety property
for a system consists in ensuring that it will never reach a dangerous
or an unwanted configuration. Safety verification is usually translated
into a reachability analysis problem \cite{frehse2011spaceex,asarin2000approximate,sun2002controllability,tiwari2003approximate,chutinan1999verification}.
Starting from an initial region, a system must not reach some unsafe
region. Different methods have been considered to address this problem.
One may explicitly compute the reachable region and determine
whether the system reaches the unsafe region \cite{gulwani2008constraint}.
An alternative idea is to compute an invariant for the system, \emph{i.e.},
a region in which the system is guaranteed to stay \cite{chutinan1999verification}.
This paper considers a class of invariants determined by \emph{barrier
functions}.

A barrier function \cite{prajna2004safety,prajna2006barrier} partitions
the state space and isolates an unsafe region from the part of the
state space containing the initial region. In \cite{prajna2004safety}
polynomial barriers are considered for polynomial systems and semi-definite
programming is used to find satisfying barrier functions. Our aim
is to extend the class of considered problems to non-polynomial systems
and to non-polynomial barriers. This paper focuses on continuous-time
systems.

The design of a barrier function is formulated as a quantified constraints
satisfaction problem (QCSP) \cite{benhamou2000universally,ratschan2006efficient}.
Interval analysis is then used to find the parameters of a barrier
function such that the QCSP is satisfied. More specifically, the algorithm
presented in \cite{Jaulin1996guaranteed} for robust controller design
is adapted and supplemented with some of the pruning schemes found
in \cite{chabert2009contractor} to solve the QCSP associated to the
barrier function design.

The paper is organized as follows. Section~\ref{sec:related} introduces
some related work. Section~\ref{sec:formulation} defines the notion
of barrier functions and formulates the design of barrier functions
as a QCSP. Section~\ref{sec:algorithm} presents the framework developed
to solve the QCSP. Design examples are presented in Section~\ref{sec:example}.
Section~\ref{sec:conclusion} concludes the work.

In what follows small italic letters $x$ represent real variables
while real vectors $\mathbf{x}$ are in bold. Intervals $\left[x\right]$
and interval vectors (boxes) $\mathbf{\left[x\right]}$ are represented
between brackets. We denote by $\mathbb{IR}$ the set of closed intervals
over $\mathbb{R}$, the set of real numbers. Data structures or sets
$\mathcal{S}$ are in upper-case calligraphic. The derivative of a
function $x$ with respect to time $t$ is denoted by~$\dot{x}$.

\section{Related work}

\label{sec:related} To prove the safety of a dynamical system, different approaches have
been proposed~\cite{Gueguen20041253,tomlin1998conflict,lygeros1999controllers,asarin2000approximate,frehse2011spaceex,tiwari2003approximate,chutinan1999verification}.
One way is to explicitly compute an outer approximation of the reachable
region from the initial region, \emph{i.e.}, the set of possible values
for the initial state. If it does not intersect the unsafe region,
then the system is safe. In \cite{asarin2000approximate,tiwari2003approximate,girard2006efficient}
the reachable region is computed for linear hybrid systems for a finite
time horizon using geometric representation such as polyhedra. The
reachable region for non-linear systems is computed in \cite{sankaranarayanan2004constructing}
using an abstraction of the non-linear systems by a linear system
expressed in a new system of coordinates. The reachability of a non-linear
system is formulated as an optimization problem in \cite{chutinan2003computational}.
In \cite{chen2012taylor}, the Picard-Lindelöf operator is combined
with Taylor models to find the reachable region for non-linear hybrid
systems. The main downside of the reachability approach is the introduction
of over-approximations during the computations which may lead to difficulties
to decide whether the system is safe.

An alternative way to address the safety problem is by exhibiting
an invariant region in which the system remains. If the invariant
does not intersect the unsafe region then the safety of the system
is proved. One way to find such an invariant is by using stability
properties of the considered dynamical system \cite{genesio1985estimation}
and to search for a Lyapunov function. In \cite{parrilo2003semidefinite}
a sum of squares decomposition and a semi-definite programming approach
are employed to find a Lyapunov function for a system with polynomial
dynamics. A template approach is considered in \cite{ratschan2006providing}
to find Lyapunov functions using a branch and relax scheme and linear
programming to solve the induced constraints. A more general idea
about invariants is introduced in \cite{prajna2004safety}. Instead
of looking for a function that fulfills some stability conditions,
a function is searched that separates the initial region from the
unsafe region. This idea is extended in \cite{gulwani2008constraint}
to search for invariants in conjunctive normal form for hybrid systems.

\section{Formulation}

\label{sec:formulation} This section recalls the safety characterization introduced in \cite{prajna2004safety}
for continuous-time systems using barrier functions.

\subsection{Safety for continuous-time systems}

Consider the autonomous continuous-time perturbed dynamical system
\begin{equation}
\mathbf{\dot{x}}=f(\mathbf{x},\mathbf{d}),\label{eq:ODE}
\end{equation}
where $\mathbf{x}\in\mathcal{X}\subseteq\mathbb{R}^{n}$ is the state
vector and $\mathbf{d}\in\mathcal{D}$ is a constant and bounded
disturbance. The set of possible initial states at $t=0$ is denoted
$\mathcal{X}_{0}\subset\mathcal{X}$. There is some unsafe subset
$\mathcal{X}_{u}\subseteq\mathcal{X}$ that shall not be reached by the
system, whatever $\mathbf{x}_{0}\mathcal{\in X}_{0}$ at time $t=0$ and
whatever $\mathbf{d}\in\mathcal{D}$. We assume that classical
hypotheses (see, \emph{e.g.,} \cite{bellman1963differential}) on $f$
are satisfied so that~\eqref{eq:ODE} has a unique solution
$\mathbf{x}(t,\mathbf{x}_{0},\mathbf{d})\in\mathcal{X}$ for any given
initial value $\mathbf{x}_{0}\mathcal{\in X}_{0}$ at time \emph{$t=0$}
and any $\mathbf{d}\in\mathcal{D}$.

\begin{definition}The dynamical system \eqref{eq:ODE} is \emph{safe}
if $\forall\mathbf{x}_{0}\in\mathcal{X}_{0}$, $\forall\mathbf{d}\in\mathcal{D}$
and $\forall t\geqslant0$, $\mathbf{x}(t,\mathbf{x}_{0},\mathbf{d})\notin\mathcal{X}_{u}$.
\end{definition}

\subsection{Barrier certificates}

A way to prove that \eqref{eq:ODE} is safe is by the barrier certificate
approach introduced in  \cite{prajna2004safety}. A barrier is a differentiable
function $B:\mathcal{X}\to\mathbb{R}$ that partitions the state space
$\mathcal{X}$ into $\mathcal{X}_{-}$ where $B(\mathbf{x})\leqslant0$
and $\mathcal{X}_{+}$ where $B(\mathbf{x})>0$ such that $\mathcal{X}_{0}\subseteq\mathcal{X}_{-}$
and $\mathcal{X}_{u}\subseteq\mathcal{X}_{+}$. Moreover, $B$ has
to be such that $\forall\mathbf{x}_{0}\in\mathcal{X}_{0}$, $\forall\mathbf{d}\in\mathcal{D}$,
$\forall t\geqslant0$, $B(\mathbf{x}(t,\mathbf{x}_{0},\mathbf{d}))\leqslant0$.

Proving that $B(\mathbf{x}(t,\mathbf{x}_{0},\mathbf{d}))\leqslant0$
requires an evaluation of the solution of~\eqref{eq:ODE} for all
$\mathbf{x}_{0}\in\mathcal{X}_{0}$ and $\mathbf{d}\in\mathcal{D}$.
Alternatively, \cite{prajna2004safety} provides some sufficient conditions
a barrier function has to satisfy to prove the safety of a dynamical
system, see Theorem~\ref{thm:Barrier_certificate}.
\begin{thm}
\label{thm:Barrier_certificate}
Consider the dynamical system \eqref{eq:ODE} and the sets
$\mathcal{X}$, $\mathcal{D}$, $\mathcal{X}_{0}$ and
$\mathcal{X}_{u}$. If there exists a function
$B:\mathcal{X}\to\mathbb{R}$ such that
\begin{equation}
  \forall\mathbf{x\in}\mathcal{X}_{0},\quad B(\mathbf{x})\leqslant 0,
  \label{eq:barrier-condition-init}
\end{equation}
\begin{equation}
  \forall\mathbf{x\in}\mathcal{X}_{u},\quad B(\mathbf{x})>0,
  \label{eq:barrier-condition-unsafe}
\end{equation}
$\forall\mathbf{x}\in\mathcal{X},\,\forall\mathbf{d\in\mathcal{D}}$,
\begin{equation}
B(\mathbf{x})=0\implies\left\langle
    \frac{\partial B(\mathbf{x})}{\partial\mathbf{x}},f(\mathbf{x},\mathbf{d})
  \right\rangle <0,
  \label{eq:barrier-condition-state}
\end{equation}
then \eqref{eq:ODE} is safe.
\end{thm}
In \eqref{eq:barrier-condition-state} $\left\langle .,.\right\rangle$
stands for the dot product in $\mathbb{R}^{n}$. In
Theorem~\ref{thm:Barrier_certificate},
\eqref{eq:barrier-condition-init} and
\eqref{eq:barrier-condition-unsafe} ensure that
$\mathcal{X}_{0}\subseteq\mathcal{X}_{-}$, and
$\mathcal{X}_{u}\subseteq\mathcal{X}_{+}$, while
\eqref{eq:barrier-condition-state} states that if $\mathbf{x}$ is on
the border between $\mathcal{X}_{-}$ and $\mathcal{X}_{+}$
(\emph{i.e.}, $B(\mathbf{x})=0$), then the dynamics $f$ pushes the
state back in $\mathcal{X}_{-}$ whatever the value of the
disturbance~$\mathbf{d}$.

\subsection{Parametric barrier functions}

The search for a barrier $B$ is challenging since it is over a
functional space. As in \cite{prajna2004safety}, this paper considers
barriers belonging to a family of parametric functions (or templates)
$B(\mathbf{x},\mathbf{p})$ depending on a parameter vector
\textbf{$\mathbf{p}\in\mathcal{P}\subseteq\mathbb{R}^{m}$}.  Then one
may search for some parameter value $\mathbf{p}$ such that
$B(\mathbf{x},\mathbf{p})$ satisfies
\eqref{eq:barrier-condition-init}-\eqref{eq:barrier-condition-state}.

If there is no \textbf{$\mathbf{p}\in\mathcal{P}$} such that
$B(\mathbf{x},\mathbf{p})$ satisfies
\eqref{eq:barrier-condition-init}-\eqref{eq:barrier-condition-state},
this does not mean that the system is not safe: other structures of
functions $B(\mathbf{x},\mathbf{p})$ could provide a barrier
certificate.

\section{Characterization using interval analysis}

\label{sec:algorithm} This section presents an approach to find a barrier function that
fulfills the constraints of Theorem~\ref{thm:Barrier_certificate}.
These constraints are first reformulated to cast the design of a
barrier function as a quantified constraint satisfaction problem
(QCSP)\cite{ratschan2006efficient}.

\subsection{Constraint satisfaction problem}

Assume that there exist some functions
$g_{0}:\mathcal{X}\to\mathbb{R}$ and $g_{u}:\mathcal{X}\to\mathbb{R}$,
such that
\begin{equation}
  \mathcal{X}_{0}=\{\mathbf{x}\in\mathcal{X}\,|\,g_{0}(\mathbf{x})\leqslant0\}
  \label{eq:X0def}
\end{equation}
and
\begin{equation}
  \mathcal{X}_{u}=\{\mathbf{x}\in\mathcal{X}\,|\,g_{u}(\mathbf{x})\leqslant0\}.
  \label{eq:Xudef}
\end{equation}

Theorem~\ref{thm:Barrier_certificate} may be reformulated as follows.
\begin{prop}
  \label{constraint-simplified}
  If $\exists\mathbf{p}\in\mathcal{P}$ such that
  $\forall\mathbf{x}\in\mathcal{X}$,
  $\forall\mathbf{d}\in\mathcal{D}$
\begin{align}
  \xi\left(\mathbf{x},\mathbf{p},\mathbf{d}\right) = &
  \left(g_{0}(\mathbf{x})>0\lor B(\mathbf{x},\mathbf{p})\leqslant0\right)
  \label{eq:TestInit1}
  \\
  &\hspace{-1.5cm} \land \left(g_{u}(\mathbf{x})>0\lor B(\mathbf{x},\mathbf{p})>0\right)
  \label{eq:TestUnsafe1}
  \\
  &\hspace{-1.5cm} \land  \left(B(\mathbf{x},\mathbf{p})\neq0\lor\left\langle
      \frac{\partial B}{\partial\mathbf{x}}(\mathbf{x},\mathbf{p}),
      f(\mathbf{x},\mathbf{d})\right\rangle <0\right)
  \label{eq:TestBorder1}
\end{align}
holds true, then the dynamical system \eqref{eq:ODE} is safe. \end{prop}
\begin{pf}
The first component of $\xi\left(\mathbf{x},\mathbf{p},\mathbf{d}\right)$,
\begin{equation}
  \xi_{\text{0}}\left(\mathbf{x},\mathbf{p}\right)=
  \left(g_{0}(\mathbf{x})>0\lor B(\mathbf{x},\mathbf{p})\leqslant0\right)
  \label{eq:TestInit}
\end{equation}
may be rewritten as
\[
\xi_{\text{0}}\left(\mathbf{x},\mathbf{p}\right)=
\left(g_{0}(\mathbf{x})\leqslant 0 \implies
  B(\mathbf{x},\mathbf{p})\leqslant0\right),
\]
see, \emph{e.g.}, \cite{gallier1986logic}. If
$\xi_{\text{0}}\left(\mathbf{x},\mathbf{p}\right)$ holds true for some
$\mathbf{p}\in\mathcal{P}$ and
$\mathbf{x}\in\mathcal{X}$, then one has either
$\mathbf{x}\in\mathcal{X}_{0}$ and
$B(\mathbf{x},\mathbf{p})\leqslant0$, or
$\mathbf{x}\notin\mathcal{X}_{0}$. In both cases,
\eqref{eq:barrier-condition-init} is satisfied. A similar derivation
can be made for the second component of
$\xi\left(\mathbf{x},\mathbf{p},\mathbf{d}\right)$ to encode
\eqref{eq:barrier-condition-unsafe},
\begin{equation}
  \xi_{u}(\mathbf{x},\mathbf{p})=
  \left(g_{u}(\mathbf{x})>0\lor B(\mathbf{x},\mathbf{p})>0\right).
  \label{eq:TestUnsafe}
\end{equation}
Now, one may rewrite the last component of
$t\left(\mathbf{x},\mathbf{p},\mathbf{d}\right)$,
\begin{equation}
  \xi_{b}(\mathbf{x},\mathbf{p},\mathbf{d)}=
  \left(B(\mathbf{x},\mathbf{p})\neq0\lor\left\langle
      \frac{\partial B}{\partial x}(\mathbf{x},\mathbf{p}),
      f(\mathbf{x},\mathbf{d})\right\rangle <0\right)
  \label{eq:TestBorder}
\end{equation}
as
\begin{multline}
  \xi_{b}(\mathbf{x},\mathbf{p},\mathbf{d)}= \\
  \left(B(\mathbf{x},\mathbf{p})=0\implies
    \left\langle
        \frac{\partial B}{\partial\vec{x}}(\mathbf{x},\mathbf{p}),
        f(\mathbf{x},\mathbf{d})\right\rangle <0\right),
  \label{eq:state-component-imply}
\end{multline}
which corresponds to \eqref{eq:barrier-condition-state}. If $\exists \mathbf{p}\in\mathcal{P}$ such that
$\forall\mathbf{x}\in\mathcal{X}$,
$\forall\mathbf{d}\in\mathcal{D}$,
$\xi\left(\mathbf{x},\mathbf{p},\mathbf{d}\right)$ holds true, then
the conditions of Theorem~\ref{thm:Barrier_certificate} are satisfied
and \eqref{eq:ODE} is safe. \hfill$\qed$

In \cite{prajna2004safety}, \eqref{eq:TestBorder1} is relaxed into
\begin{equation}
  \xi_{\text{b}}\left(\mathbf{x},\mathbf{p},\mathbf{d}\right)=
  \left(\left\langle \frac{\partial B}{\partial\vec{x}}(\mathbf{x},\mathbf{p}),
      f(\mathbf{x},\mathbf{d})\right\rangle <0\right),
  \label{eq:Tborderrelaxed}
\end{equation}
with the consequence of possible elimination of barrier functions that
would satisfy \eqref{eq:TestBorder1} for some $\mathbf{p}$ but not
\eqref{eq:Tborderrelaxed}. Our aim in this paper is to design barrier
functions without resorting to this relaxation by considering methods
from interval analysis \cite{jaulin2001applied} which allow to
consider strongly nonlinear dynamics and barrier functions.
\end{pf}

\subsection{Solving the constraints}

\label{sub:SolvingConstraints}

To find a valid barrier function one needs to find some
$\mathbf{p}\in\mathcal{P}$ such that
$B(\mathbf{x},\mathbf{p})$ satisfies the conditions of
Proposition~\ref{constraint-simplified}. For that purpose, the
\emph{Computable Sufficient Conditions}-\emph{Feasible Point Searcher}
(CSC-FPS) algorithm \cite{Jaulin1996guaranteed} is adapted.

In what follows, we assume that $\mathcal{X}$, $\mathcal{D}$, and $\mathcal{P}$ are boxes, \emph{i.e.},
$\mathcal{X}=\left[\mathbf{x}\right]$, $\mathcal{D}=\left[\mathbf{d}\right]$, and
$\mathcal{P}=\left[\mathbf{p}\right]$. CSC-FPS may also be applied when $\mathcal{X}$, $\mathcal{D}$, and $\mathcal{P}$
consist of a union of non-overlapping boxes.

Consider some function $g:\mathbb{R}^{n}\times\mathbb{R}^{m}\rightarrow\mathbb{R}^{k}$ and
some box $\left[\mathbf{z}\right]\in\mathbb{IR}^{k}$. CSC-FPS is
designed to determine whether
\begin{equation}
  \exists\mathbf{p}\in\left[\mathbf{p}\right],
  \mbox{ }\forall\mathbf{x}\in\left[\mathbf{x}\right],
  \mbox{ }g\left(\mathbf{x},\mathbf{p}\right)\in
  \left[\mathbf{z}\right]
  \label{eq:GenercRobustControl}
\end{equation}
and to provide some satisfying $\mathbf{p}$. We extend CSC-FPS to
handle conjunctions and disjunctions of constraints and supplement it
with efficient pruning techniques involving contractors provided by
interval analysis \cite{jaulin2001applied}.

FPS branches over the parameter search box $\left[\mathbf{p}\right]$.
Branching is performed based on the results provided by CSC. For a
given box
$\mathbf{\left[p\right]}_{0}\subseteq\left[\mathbf{p}\right]$, CSC
returns \texttt{true} when it manages to prove that
\eqref{eq:GenercRobustControl} is satisfied for some
$\mathbf{p}\in\mathbf{\left[p\right]}_{0}$.  CSC returns
\texttt{false} when it is able to show that there is no
$\mathbf{p}\in\mathbf{\left[p\right]}_{0}$ satisfying
\eqref{eq:GenercRobustControl}. CSC returns \texttt{unknown} in the
other cases.

In Proposition~\ref{constraint-simplified},
$\xi\left(\mathbf{x},\mathbf{p},\mathbf{d}\right)$ consists of the
conjunction of three terms of the form
\begin{equation}
  \tau\left(\mathbf{x},\mathbf{p},\mathbf{d}\right)=
  \left(u(\mathbf{x},\mathbf{p})\in\mathcal{A}\right)\lor
  \left(v(\mathbf{x},\mathbf{p,d})\in\mathcal{B}\right).
  \label{eq:ElementaryTest}
\end{equation}
For $\xi_{0}(\mathbf{x},\mathbf{p})$, defined in \eqref{eq:TestInit1},
\begin{equation}
  \mathcal{A}=\left]0,+\infty\right]\mbox{ and }
  \mathcal{B}=\left]-\infty,0\right];
  \label{eq:ABInit}
\end{equation}
for $\xi_{u}(\mathbf{x},\mathbf{p})$, defined in
\eqref{eq:TestUnsafe1},
\begin{equation}
  \mathcal{A}=\left]0,+\infty\right]\mbox{ and }
  \mathcal{B}=\left]0,+\infty\right[;
  \label{eq:ABUnsafe}
\end{equation}
for $\xi_{b}(\mathbf{x},\mathbf{p},\mathbf{d})$, defined in
\eqref{eq:TestBorder1},
\begin{equation}
  \mathcal{A}=\left]-\infty,0\right[\cup\left]0,+\infty\right[
  \mbox{ and }
  \mathcal{B}=\left]-\infty,0\right[.
  \label{eq:ABBorder}
\end{equation}

To illustrate the main ideas of CSC-FPS combined with contractors, one
focuses on the generic QCSP
\begin{equation}
  \exists\mathbf{p}\in\left[\mathbf{p}\right],
  \forall\mathbf{x}\in\mathbf{\left[x\right],
    \forall d\in}\left[\mathbf{d}\right],\,
  \tau\left(\mathbf{x},\mathbf{p},\mathbf{d}\right)\mbox{ holds true}.
  \label{eq:tau-QSCP}
\end{equation}
Finding a solution for such QCSP involves three steps: validation,
reduction of the parameter and state spaces, and bisection.

\subsubsection{Validation}
\label{sub:Validation}

In the validation step, one tries to prove that some vector
$\mathbf{p}\in\mathbf{\left[p\right]}$ is such that
$\forall\mathbf{x}\in\mathbf{\left[x\right]}$,
$\forall\mathbf{d\in}\left[\mathbf{d}\right]$,
$\tau\left(\mathbf{x},\mathbf{p,d}\right)$ holds true. By definition
of $\tau\left(\mathbf{x},\mathbf{p},\mathbf{d}\right)$, one has to
prove that
\begin{multline}
  \exists\mathbf{p}\in\left[\mathbf{p}\right],
  \forall\mathbf{x}\in\left[\mathbf{x}\right],
  \forall\mathbf{d}\in\left[\mathbf{d}\right],\\
  \left(u(\mathbf{x},\mathbf{p})\in\mathcal{A}\right) \lor
  \left(v(\mathbf{x},\mathbf{p},\mathbf{d})\in\mathcal{B}\right).
  \label{eq:Validation}
\end{multline}
For that purpose, one chooses some arbitrary
$\mathbf{p}\in\mathbf{\left[p\right]}$ and evaluates the set of values
$u\left(\left[\mathbf{x}\right],\mathbf{p}\right)=\left\{
  u\left(\mathbf{x},\mathbf{p}\right)\mid\mathbf{x}\in\mathbf{\left[\mathbf{x}\right]}\right\}
$ and
$v\left(\left[\mathbf{x}\right],\mathbf{p,\left[d\right]}\right)=\{
v\left(\mathbf{x},\mathbf{p},\mathbf{d}\right)\mid\mathbf{x}\in\mathbf{\left[\mathbf{x}\right]},\mathbf{d}\in\mathbf{\left[d\right]}\}
$ taken by $u(\mathbf{x},\mathbf{p})$ and $v(\mathbf{x},\mathbf{p,d})$
for all $\mathbf{x}\in\mathbf{\left[\mathbf{x}\right]}$ and for all
$\mathbf{d}\in\left[\mathbf{d}\right]$. Outer-approximations of
$u\left(\left[\mathbf{x}\right],\mathbf{p}\right)$ and
$v\left(\left[\mathbf{x}\right],\mathbf{p,\left[d\right]}\right)$ are
easily obtained using \emph{inclusion functions} provided by interval
analysis.
\begin{defn}
  An inclusion function $[f]:\mathbb{IR}^{n}\to\mathbb{IR}^{k}$ for a
  function $f:\mathbb{R}^{n}\to\mathbb{R}^{k}$ satisfies for all
  $\left[\mathbf{x}\right]\in\mathbb{IR}^{n}$,
  \begin{equation}
    f\left(\left[\mathbf{x}\right]\right) =
    \left\{ f(\mathbf{x})\mid\mathbf{x}\in\left[\mathbf{x}\right]\right\}
    \subseteq\left[f\right]\left(\left[\mathbf{x}\right]\right).
    \label{eq:Inclusion_function}
  \end{equation}
\end{defn}

The \emph{natural} inclusion function is the simplest to obtain: all
occurrences of the real variables are replaced by their interval
counterpart and all arithmetic operations are evaluated using interval
arithmetic. More sophisticated inclusion functions such as the
centered form, or the Taylor inclusion function may also be used, see
\cite{jaulin2001applied}.

Using inclusion functions, one may evaluate whether
\begin{displaymath}
  [u](\mathbf{\left[x\right]},\mathbf{p})\subseteq\mathcal{A}
  \text{\quad or\quad}
  [v](\mathbf{\left[x\right]},\mathbf{p,\left[d\right]})\subseteq\mathcal{B}
  \quad\text{holds true}
\end{displaymath}
for the various $\mathcal{A}$ and $\mathcal{B}$ defined in
\eqref{eq:ABInit}, \eqref{eq:ABUnsafe}, and \eqref{eq:ABBorder}.

Different choices can be considered for $\mathbf{p}$: one can take a
random point in $\mathbf{\left[p\right]}$, the middle, or one of the
edges of $\mathbf{\left[p\right]}$. Here, only the middle of
$\mathbf{\left[p\right]}$ is considered.

\subsubsection{Reduction of the parameter and state spaces}

\label{sub:Reduction}

To facilitate the search for $\mathbf{p}\in\mathbf{\left[p\right]}$
one may previously eliminate parts of $\mathbf{\left[p\right]}$ which
may be proved not to contain any $\mathbf{p}$ satisfying
\eqref{eq:Validation}. The elimination process can be done by
evaluation or by using
\emph{contractors}~\cite{chabert2009contractor}.

\paragraph{Evaluation}

From \eqref{eq:Validation} one deduces that a box
$\mathbf{\left[p\right]}$ can be eliminated, \emph{i.e.}, shown not to
contain any $\mathbf{p}$ satisfying \eqref{eq:ElementaryTest}, if
\begin{multline}
  \forall\mathbf{p}\in\mathbf{\left[p\right]},\exists
  \mathbf{x}\in\left[\mathbf{x}\right],
  \exists\mathbf{d}\in\left[\mathbf{d}\right],\\
  u(\mathbf{x},\mathbf{p})\subseteq\overline{\mathcal{A}} \wedge
  v(\mathbf{x},\mathbf{p,d})\subseteq\overline{\mathcal{B}},
  \label{eq:invalidation}
\end{multline}
where $\mathcal{\overline{A}}=\mathbb{R}\setminus\mathcal{A}$ and
$\mathcal{\overline{B}}=\mathbb{R}\setminus\mathcal{B}$. Using again
an inclusion function, one may verify whether
$[u](\mathbf{x},\mathbf{\left[p\right]})\subseteq\overline{\mathcal{A}}$
and
$[v](\mathbf{x},\mathbf{\left[p\right],d})\subseteq\overline{\mathcal{B}}$
for some $\mathbf{x}\in\mathbf{\left[x\right]}$ and
$\mathbf{d\in\left[d\right]}$ and thus eliminate the box
$\mathbf{\left[p\right]}$.

One may easily verify \eqref{eq:TestInit1} and \eqref{eq:TestUnsafe1}
since $\overline{\mathcal{A}}$ and $\overline{\mathcal{B}}$ are
unbounded and may contain $[u](\mathbf{x},\mathbf{\left[p\right]})$
and $[v](\mathbf{x},\mathbf{\left[p\right],d})$. For
\eqref{eq:TestBorder1} the inclusion is impossible to prove except in
degenerate cases since $\overline{\mathcal{A}}=\left\{ 0\right\} $.

\paragraph{Contractors}

\label{par:contractor}

Consider some function $g:\mathbb{R}^{n}\rightarrow\mathbb{R}^{k}$
and some set $\mathcal{Z}\subset\mathbb{R}^{k}$.
\begin{defn}
  A contractor $\mathcal{C}_{c}:\mathbb{IR}^{n}\to\mathbb{IR}^{n}$
  associated to the generic constraint
  \begin{equation}
    c:g(\mathbf{x})\in\mathcal{Z}\label{eq:elementary-constraint}
  \end{equation}
  is a function taking a box $\left[\mathbf{x}\right]$ as input and
  returning a box
  $\mathcal{C}_{c}\left(\left[\mathbf{x}\right]\right)$ satisfying
  \begin{equation}
    \mathcal{C}_{c}\left(\left[\mathbf{x}\right]\right)
    \subseteq\left[\mathbf{x}\right]
    \label{eq:Contractor_Contraction}
  \end{equation}
  and
  \begin{equation}
    g\left(\left[\mathbf{x}\right]\right)\cap\mathcal{Z}=
    g\left(\mathcal{C}_{c}\left(\left[\mathbf{x}\right]\right)\right)
    \cap\mathcal{Z}.
    \label{eq:Contractor_NoLoss}
  \end{equation}
\end{defn}
$\mathcal{C}_{c}$ provides a box containing the solutions of
$g(\mathbf{x})\in\mathcal{Z}$ included in $\left[\mathbf{x}\right]$:
\eqref{eq:Contractor_Contraction} ensures that the returned box is
included in $\mathbf{\left[x\right]}$ and \eqref{eq:Contractor_NoLoss}
ensures that no solution of $g(\mathbf{x})\in\mathcal{Z}$ in
$\mathbf{\left[x\right]}$ is lost.

Consider now two functions
$g_{1}:\mathbb{R}^{n}\rightarrow\mathbb{R}^{k_{1}}$ and
$g_{2}:\mathbb{R}^{n}\rightarrow\mathbb{R}^{k_{2}}$, two sets
$\mathcal{Z}_{1}\subset\mathbb{R}^{k_{1}}$ and
$\mathcal{Z}_{2}\subset\mathbb{R}^{k_{2}}$, and the associated
constraints
$c_{1}:g_{1}(\mathbf{\left[x\right]})\subseteq\mathcal{Z}_{1}$ and
$c_{2}:g_{2}(\mathbf{\left[x\right]})\subseteq\mathcal{Z}_{2}$.
Assume that two contractors $\mathcal{C}_{c_{1}}$ and
$\mathcal{C}_{c_{2}}$ are available for $c_{1}$ and $c_{2}$. A
contractor $\mathcal{C}_{c_{1}\land c_{2}}$ for the conjunction
$c_{1}\wedge c_{2}$ of $c_{1}$ and $c_{2}$ may be obtained as
\begin{equation}
  \mathcal{C}_{c_{1}\land c_{2}}(\mathbf{\left[x\right]})=
  \mathcal{C}_{c_{1}}(\mathbf{\left[x\right]})\cap
  \mathcal{C}_{c_{2}}(\mathbf{\left[x\right]}),
  \label{eq:conj_contracteur}
\end{equation}
or by composition of contractors
\begin{equation}
  \mathcal{C}_{c_{1}\land c_{2}}(\mathbf{\left[x\right]})=
  \mathcal{C}_{c_{2}}(\mathcal{C}_{c_{1}}(\mathbf{\left[x\right]})).
  \label{eq:conj_contracteur2}
\end{equation}
A contractor $\mathcal{C}_{c_{1}\vee c_{2}}$ for the disjunction
$c_{1}\vee c_{2}$ of $c_{1}$ and $c_{2}$ may be obtained as follows
\begin{equation}
  \mathcal{C}_{c_{1}\lor c_{2}}(\mathbf{\left[x\right]})=
  \Box\{\mathcal{C}_{c_{1}}(\mathbf{\left[x\right]})\cup\mathcal{C}_{c_{2}}
  (\mathbf{\left[x\right]})\},
  \label{eq:disj_contracteur}
\end{equation}
see \cite{chabert2009contractor}, with $\Box\{\cdot\}$ the interval
hull of a set.

\begin{figure}
\centering
\includegraphics[width=1\columnwidth]{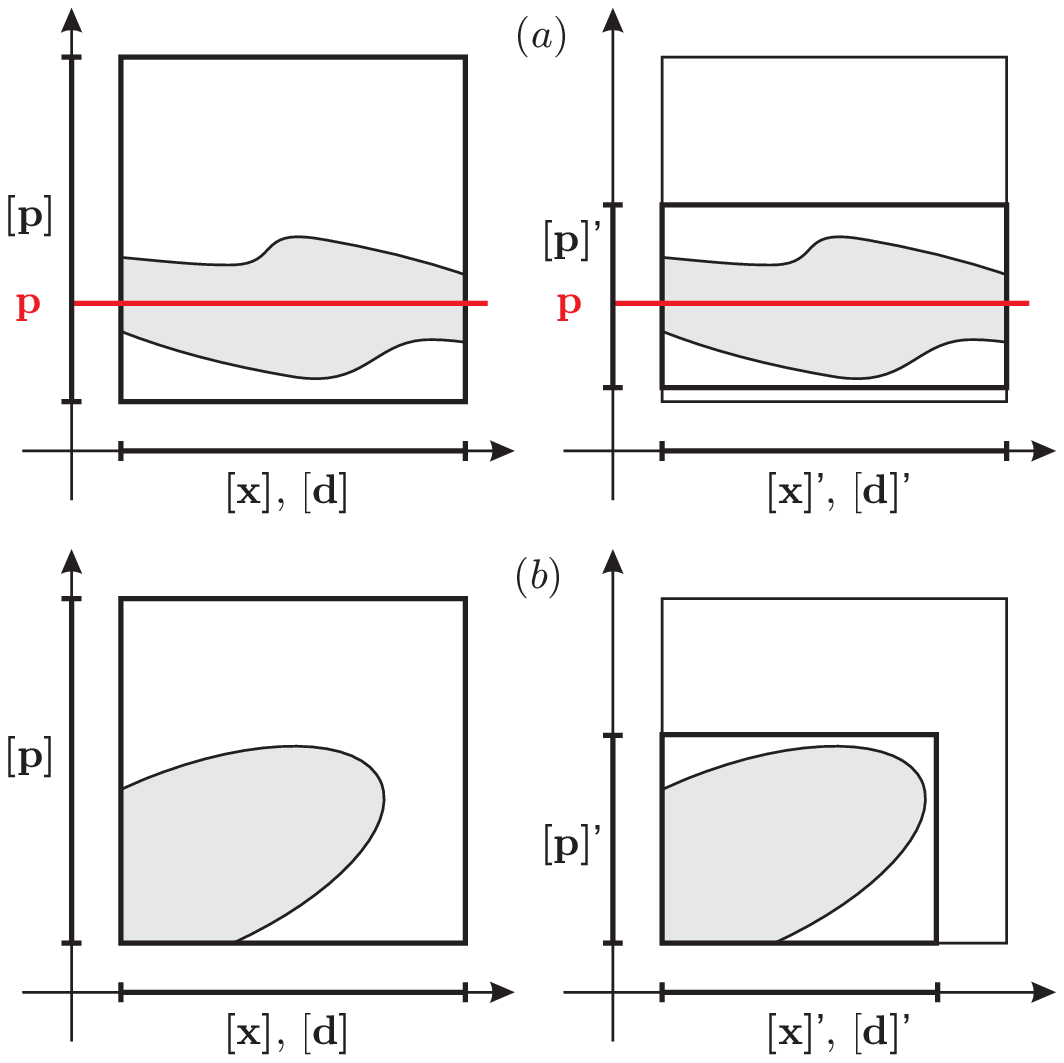}\protect
\caption{Contractions using $\mathcal{C}_{\tau}$; the set for which \eqref{eq:ElementaryReducing}
is satisfied is in gray; (\emph{a}) $\mathbf{\left[p\right]}\setminus\mathbf{\left[p\right]}'\protect\neq\emptyset$
and the search space for satisfying $\mathbf{p}$ can be reduced to
$\left[\mathbf{p}\right]';$ (\emph{b}) $\mathbf{\left[x\right]}'\protect\neq\mathbf{\left[x\right]}$
or $\mathbf{\left[d\right]}'\protect\neq\mathbf{\left[d\right]}$,
it is thus not possible to find some \textbf{$\mathbf{p}\in\left[\mathbf{p}\right]$}
such that \eqref{eq:ElementaryReducing} is satisfied for all\textbf{
$\mathbf{x}\in\left[\mathbf{x}\right]$} and all $\mathbf{d}\in\left[\mathbf{d}\right]$.
\label{fig:Contractions}}
\end{figure}

Using a contractor $\mathcal{C}_{c}$ for
\eqref{eq:elementary-constraint}, one is able to characterize some
$\left[\tilde{\mathbf{x}}\right]\subset\left[\mathbf{x}\right]$ such
that $\forall\mathbf{x}\in\left[\tilde{\mathbf{x}}\right]$,
$g(\mathbf{x})\notin\mathcal{Z}$.
\begin{prop}
  \label{prop:ContractCompl}Consider a box $\mathbf{\left[x\right]}$,
  the elementary constraint \eqref{eq:elementary-constraint}, and the
  contracted box $\mathcal{C}_{c}\left(\mathbf{\left[x\right]}\right)\subseteq\left[\mathbf{x}\right]$.
  Then,
  \begin{equation}
    \forall\mathbf{x}\in\left[\mathbf{x}\right]\backslash\mathcal{C}_{c}\left(\mathbf{\left[x\right]}\right),\mbox{ one has }g(\mathbf{x})\notin\mathcal{Z},\label{eq:Contradiction}
  \end{equation}
  where $\left[\mathbf{x}\right]\backslash\mathcal{C}_{c}\left(\mathbf{\left[x\right]}\right)$
  denotes the box $\left[\mathbf{x}\right]$ deprived from $\mathcal{C}_{c}\left(\mathbf{\left[x\right]}\right)$,
  which is not necessarily a box.\end{prop}
\begin{pf}
  Consider
  $\mathbf{x}\in\left[\mathbf{x}\right]\backslash\mathcal{C}_{c}\left(\mathbf{\left[x\right]}\right)$
  and assume that $g(\mathbf{x})\in\mathcal{Z}$. Since
  $g(\mathbf{x})\in\mathcal{Z}$ and
  $\mathbf{x}\in\left[\mathbf{x}\right]$, one should have
  $\mathbf{x}\in\mathcal{C}_{c}\left(\mathbf{\left[x\right]}\right)$,
  according to \eqref{eq:Contractor_NoLoss}, which contradicts the
  fact that
  $\mathbf{x}\in\left[\mathbf{x}\right]\backslash\mathcal{C}_{c}\left(\mathbf{\left[x\right]}\right)$.\hfill
  $\qed$

Proposition~\ref{prop:ContractCompl} can be used to eliminate
$\mathbf{\left[p\right]}$ or a part of $\mathbf{\left[p\right]}$ for
which it is not possible to find any $\mathbf{p}$ satisfying
\eqref{eq:Validation}. Consider the constraint
\begin{equation}
\tau:\left(u(\mathbf{x},\mathbf{p})\in\mathcal{A}\right)\lor\left(v(\mathbf{x},\mathbf{p},\mathbf{d})\in\mathcal{B}\right)\label{eq:ElementaryReducing}
\end{equation}
and a contractor $\mathcal{C}_{\tau}$ for this constraint. It involves
elementary contractors for the components of the disjunction in \eqref{eq:ElementaryReducing},
combined as in \eqref{eq:disj_contracteur}. For the boxes $\left[\mathbf{x}\right]$,
$\left[\mathbf{p}\right]$, and $\mathbf{\left[d\right]}$, one gets
\begin{equation}
  \begin{pmatrix}\left[\mathbf{x}\right]',&
    \left[\mathbf{p}\right]',&
    \mathbf{\left[d\right]}'
  \end{pmatrix}=\mathcal{C}_{\tau}\begin{pmatrix}\left[\mathbf{x}\right],&
    \left[\mathbf{p}\right],&
    \mathbf{\left[d\right]}
  \end{pmatrix},\label{eq:Contraction-tau}
\end{equation}
where $\left[\mathbf{x}\right]'$, $\left[\mathbf{p}\right]'$, and
$\left[\mathbf{d}\right]'$ are the contracted boxes. Three cases
may then be considered.\end{pf}
\begin{enumerate}
\item If
  $\mathbf{\left[p\right]}\setminus\mathbf{\left[p\right]}'\neq\emptyset$,
  then $\forall\mathbf{p}\in\mathbf{\left[p\right]\setminus}\mathbf{\left[p\right]}',\,\forall\mathbf{x}\in\mathbf{\left[x\right]},\,\forall\mathbf{d}\in\mathbf{\left[d\right]}$,
  \begin{equation}
    u(\mathbf{x},\mathbf{p})\notin\mathcal{A}\land v(\mathbf{x},\mathbf{p},\mathbf{d})\notin\mathcal{B},\label{eq:contract-params}
  \end{equation}
  and there is no
  $\mathbf{p}\in\left[\mathbf{p}\right]\setminus\mathbf{\left[p\right]}'$
  such that \eqref{eq:ElementaryReducing} holds true for all
  $\mathbf{x}\in\mathbf{\left[x\right]}$ and for all
  $\mathbf{d}\in\mathbf{\left[d\right]}$. Consequently, the search
  space for $\mathbf{p}$ can be reduced to $\mathbf{\left[p\right]}'$,
  see Figure~\ref{fig:Contractions}~(a).
\item If $\mathbf{\left[x\right]}\setminus\mathbf{\left[x\right]}'\neq\emptyset$
  then, from Proposition~\ref{prop:ContractCompl}, one has $  \forall\mathbf{p\in}\mathbf{\left[p\right]},\mathbf{\,\forall x\in}\mathbf{\left[x\right]\setminus}\mathbf{\left[x\right]}',\mathbf{\,\forall d\in}\mathbf{\left[d\right]}$,
\begin{equation}
  u(\mathbf{x},\mathbf{p})\notin\mathcal{A}\land
  v(\mathbf{x},\mathbf{p},\mathbf{d})\notin\mathcal{B},\label{eq:contract-states}
\end{equation}
and there is no $\mathbf{p}\in\left[\mathbf{p}\right]$ such that
\eqref{eq:ElementaryReducing} holds true for all
$\textbf{x}\in\left[\mathbf{x}\right]$, see
Figure~\ref{fig:Contractions} (b).
\item If
  $\mathbf{\left[d\right]}\setminus\mathbf{\left[d\right]}'\neq\emptyset$,
  then $    \forall\mathbf{p}\in\mathbf{\left[p\right]},\,\forall\mathbf{x}\in\mathbf{\left[x\right]},\,\forall\mathbf{d}\in\mathbf{\left[d\right]}\setminus\mathbf{\left[d\right]}'$,
  \begin{equation}
    u(\mathbf{x},\mathbf{p})\notin\mathcal{A}\land
    v(\mathbf{x},\mathbf{p},\mathbf{d})\notin\mathcal{B},\label{eq:contract-disturb}
\end{equation}
and there is no $\mathbf{p}\in\left[\mathbf{p}\right]$ such that
\eqref{eq:ElementaryReducing} holds true for all $\textbf{d}\in\left[\mathbf{d}\right]$,
see Figure~\ref{fig:Contractions} (b).
\end{enumerate}

\begin{figure}
\centering{}\includegraphics[width=1\columnwidth]{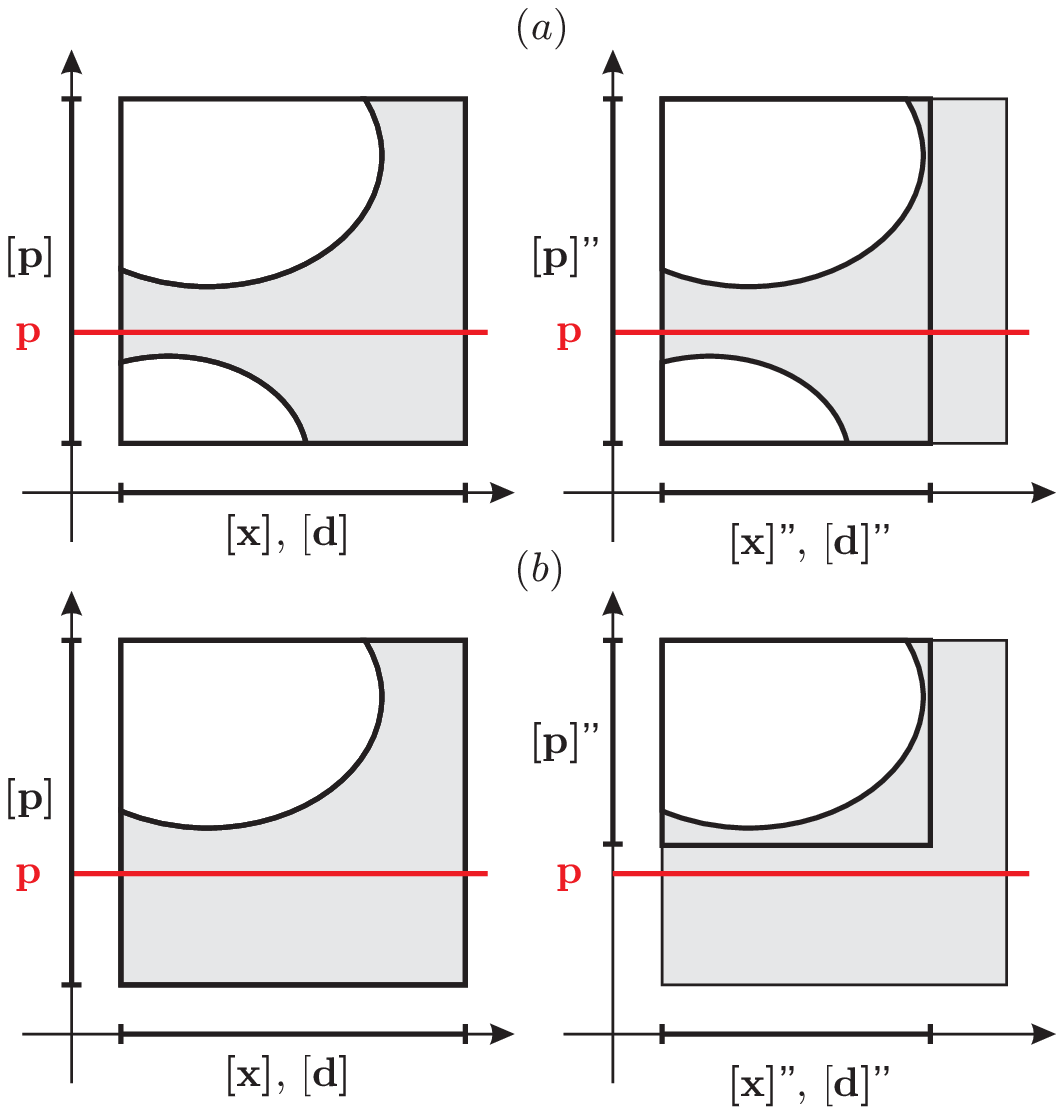}\protect\caption{Contractions using $\mathcal{C}_{\overline{\tau}}$; the set for which
\eqref{eq:Neg-ElementaryReducing} is satisfied is in white; (\emph{a})
$\mathbf{\left[x\right]}"\protect\neq\mathbf{\left[x\right]}$ and/or
$\mathbf{\left[d\right]}"\protect\neq\mathbf{\left[d\right]}$ and
one has only to find some suitable $\mathbf{p\in\left[\mathbf{p}\right]}$
such that \eqref{eq:ElementaryReducing} is satisfied for all\textbf{
$\mathbf{x}\in\left[\mathbf{x}\right]"$} and all $\mathbf{d}\in\left[\mathbf{d}\right]"$;
(\emph{b}) $\mathbf{\left[p\right]}"\protect\neq\mathbf{\left[p\right]}$,
one may choose any $\mathbf{p}\in\left[\mathbf{p}\right]\setminus\left[\mathbf{p}\right]"$
(for example the value of $\mathbf{p}$ indicated in red) and \eqref{eq:ElementaryReducing}
will hold true for all $\left(\mathbf{x},\mathbf{d}\right)\in\left(\mathbf{\left[x\right]}\times\left[\mathbf{d}\right]\right)$.
\label{fig:Contractions2}}
\end{figure}

One can reduce the size of the sets for the state $\mathbf{x}$ and
the disturbance $\mathbf{d}$ on which \eqref{eq:ElementaryReducing}
has to be verified using the contraction on the negation of this constraint.
Consider the negation $\overline{\tau}$ of $\tau$
\begin{equation}
\overline{\tau}=\left(u(\mathbf{x},\mathbf{p})\in\mathcal{\overline{A}}\right)\land\left(v(\mathbf{x},\mathbf{p},\mathbf{d})\in\mathcal{\overline{B}}\right)\label{eq:Neg-ElementaryReducing}
\end{equation}
and a contractor $\mathcal{C}_{\overline{\tau}}$ for this constraint.
Assume that after applying $\mathcal{C}_{\overline{\tau}}$ for the
boxes $\left[\mathbf{x}\right]$, $\mathbf{\left[p\right]}$, and
$\mathbf{\left[d\right]}$, one gets
\begin{equation}
\begin{pmatrix}\left[\mathbf{x}\right]",&
\left[\mathbf{p}\right]",&
\mathbf{\left[d\right]}"
\end{pmatrix}=\mathcal{C}_{\overline{\tau}}\begin{pmatrix}\left[\mathbf{x}\right],&
\left[\mathbf{p}\right],&
\mathbf{\left[d\right]}
\end{pmatrix}.\label{eq:contraction-state-space}
\end{equation}
From Proposition~\ref{prop:ContractCompl}, one knows that
\begin{multline}
  \forall\left(\mathbf{x},\mathbf{p},\mathbf{d}\right)\in\left(\mathbf{\left[x\right]}\times\mathbf{\left[p\right]}\times\left[\mathbf{d}\right]\right)\setminus\left(\left[\mathbf{x}\right]"\times\left[\mathbf{p}\right]"\times\left[\mathbf{d}\right]"\right),\\
  u(\mathbf{x},\mathbf{p})\in\mathcal{A}\lor
  v(\mathbf{x},\mathbf{p},\mathbf{d})\in\mathcal{B}.
  \label{eq:contract-neg}
\end{multline}
Indeed, if $\left[\mathbf{p}\right]"=\left[\mathbf{p}\right]$, one can
focus on the search for some $\mathbf{p}\in\mathbf{\left[p\right]}$
satisfying \eqref{eq:ElementaryReducing} by considering only
$\mathbf{\left[x\right]}"\times\left[\mathbf{d}\right]"$, since for
all
$\left(\mathbf{x},\mathbf{d}\right)\in\left(\mathbf{\left[x\right]}\times\left[\mathbf{d}\right]\right)\setminus\left(\left[\mathbf{x}\right]"\times\left[\mathbf{d}\right]"\right)$,
\eqref{eq:ElementaryReducing} is satisfied for all
$\mathbf{p}\in\mathbf{\left[p\right]}$, see
Figure~\ref{fig:Contractions2}~(\emph{a}).

Now, if $\left[\mathbf{p}\right]"\neq\left[\mathbf{p}\right]$, then
any
$\mathbf{p}\in\left[\mathbf{p}\right]\setminus\left[\mathbf{p}\right]"$
will satisfy \eqref{eq:ElementaryReducing} for all
$\left(\mathbf{x},\mathbf{d}\right)\in\left(\mathbf{\left[x\right]}\times\left[\mathbf{d}\right]\right)$,
see Figure~\ref{fig:Contractions2}~(\emph{b}).

\subsubsection{Bisection}

One is unable to decide whether some $\mathbf{p}\in\left[\mathbf{p}\right]$
satisfies \eqref{eq:Validation} when
\begin{equation}
[u](\mathbf{\left[x\right]},\mathbf{p})\cap\mathcal{A}\neq\emptyset\text{ and \ensuremath{\left[u\right]}\ensuremath{\left(\mathbf{\left[x\right]},\mathbf{p}\right)\nsubseteq\mathcal{A}}}\label{eq:DifficultyValidation}
\end{equation}
or when
\begin{equation}
[v](\mathbf{\left[x\right]},\mathbf{p,\left[d\right]})\cap\mathcal{B}\neq\emptyset\,\text{and}\,[v](\mathbf{\left[x\right]},\mathbf{p,\left[d\right]})\nsubseteq\mathcal{B}.\label{eq:DifficultValidationBis}
\end{equation}
This situation occurs in two cases. First, when the selected
$\mathbf{p}$ does not satisfy \eqref{eq:Validation} for all
$\mathbf{x}\in\left[\mathbf{x}\right]$ and for all
$\mathbf{d}\in\left[\mathbf{d}\right]$. Second, when inclusion
functions introduce some \emph{pessimism}, \emph{i.e.}, they provide
an over-approximation of the range of functions over intervals.

To address both cases, one may perform bisections of
$\left[\mathbf{x}\right]\times\left[\mathbf{d}\right]$ and try to
verify \eqref{eq:Validation} on the resulting sub-boxes for the
\emph{same} $\mathbf{p}$. Bisection allows to isolate subsets of
$\left[\mathbf{x}\right]\times\left[\mathbf{d}\right]$ on which one
may show that \eqref{eq:invalidation} holds true. Bisections also
reduce pessimism, and may thus facilitate the verification of
\eqref{eq:Validation}. The bisection of
$\left[\mathbf{x}\right]\times\left[\mathbf{d}\right]$ is performed
within CSC as long as the width of the bisected boxes are larger than
some $\varepsilon_{\text{x}}>0$. When CSC is unable to prove
\eqref{eq:Validation} or \eqref{eq:invalidation} and when all bisected
boxes are smaller than $\varepsilon_{\text{x}}$, CSC returns
\texttt{unknown}.

FPS performs similar bisections on $\mathbf{\left[\mathbf{p}\right]}$
and stops when the width of all bisected boxes are smaller than
$\varepsilon_{\text{p}}>0$.

\subsubsection{Composition of constraints}

To prove the safety of the dynamical system \eqref{eq:ODE},
Proposition~\ref{constraint-simplified} shows that one has to find
some $\mathbf{p}\in\mathbf{\left[p\right]}$ such that
$\forall\mathbf{x}\in\mathbf{\left[x\right]}$,
$\forall\mathbf{d}\in\mathbf{\left[d\right]}$,
$\xi\left(\mathbf{x},\mathbf{p},\mathbf{d}\right)$ holds true. Since
$\xi\left(\mathbf{x},\mathbf{p},\mathbf{d}\right)$ consists of the
conjunction of three elementary constraints of the form
\eqref{eq:tau-QSCP}, validation requires the verification of
\eqref{eq:Validation} for the \emph{same} \textbf{$\mathbf{p}$}
considering \eqref{eq:TestInit1}, \eqref{eq:TestUnsafe1}, and
\eqref{eq:TestBorder1} \emph{simultaneousl}y.  Invalidation may be
performed as soon as one is able to prove that one of the constraints
\eqref{eq:TestInit1}, \eqref{eq:TestUnsafe1}, or
\eqref{eq:TestBorder1} does not hold true using
\eqref{eq:invalidation}.  Contraction may benefit from the conjunction
or disjunction of these constraints, as introduced in
\eqref{eq:conj_contracteur2} and \eqref{eq:disj_contracteur}.

\subsection{CSC-FPS algorithms with contractors}

The CSC-FPS algorithm, presented in \cite{Jaulin1996guaranteed} is
supplemented with the contractors introduced in
Section~\ref{sub:SolvingConstraints} to improve its efficiency.

FPS, described in Algorithm~\ref{alg:fps}, searches for some
\emph{satisfying} $\mathbf{p}\in\left[\mathbf{p}\right]$, \emph{i.e.},
some $\mathbf{p}\in\left[\mathbf{p}\right]$ such that
$\xi\left(\mathbf{x},\mathbf{p},\mathbf{d}\right)$ introduced in
Proposition~\ref{constraint-simplified} holds true for all
$\mathbf{x}\in\mathbf{\left[x\right]}$ and all
$\mathbf{d}\in\mathbf{\left[d\right]}$. This may require to bisect
$\left[\mathbf{p}\right]$ into subboxes stored in a queue
$\mathcal{Q}$, which initial content is $\left[\mathbf{p}\right]$. A
sub-box $\left[\mathbf{p}\right]_{0}\subseteq\left[\mathbf{p}\right]$
is extracted from $\mathcal{Q}$. A reduction of
$\left[\mathbf{p}\right]_{0}$ is then performed at
Line~\ref{Alg1:Reduction} to eliminate values of
$\mathbf{p}\in\left[\mathbf{p}\right]_{0}$ which cannot be satisfying.
To facilitate contraction, specific
$\mathbf{x}\in\left[\mathbf{x}\right]$ are chosen; here only the
midpoint $\text{m}\left(\left[\mathbf{x}\right]\right)$ of
$\left[\mathbf{x}\right]$ is considered. If
$\mathbf{\left[p\right]}_{0}'$ is empty, the next box in $\mathcal{Q}$
has to be explored. Then CSC is called for each constraint $t_{0}$,
$t_{u}$, and $t_{b}$ to verify whether
$\text{m}\left(\left[\mathbf{p}\right]_{0}'\right)$, the midpoint of $\left[\mathbf{p}\right]_{0}'$, is
satisfying. When all calls of CSC return \texttt{true} at
Line~\ref{alg1:Alltrue}, a barrier function with parameter
$\text{m}\left(\left[\mathbf{p}\right]_{0}'\right)$ is found. When one
of the calls of CSC returns \texttt{false} at
Line~\ref{Alg1:OneFalse}, $\left[\mathbf{p}\right]_{0}'$ is proved not
to contain any satisfying $\mathbf{p}$. In all other cases, when $\text{w}(\mathbf{\left[p\right]}_{0}')$, the
width of $\mathbf{\left[p\right]}_{0}'$, is larger than
$\varepsilon_{\text{p}}$, $\mathbf{\left[p\right]}_{0}'$ is bisected
and the resulting subboxes are stored in $\mathcal{Q}$. When
$\mathbf{\left[p\right]}_{0}'$ is too small, even if one was not able
to decide whether it constains a satisfying $\mathbf{p}$, it is not
further considered to ensure termination of FPS in finite time. The
price to be paid in such situations is the impossibility to conclude
that the initial box $\left[\mathbf{p}\right]$ does not contain some
satisfying $\mathbf{p}$. This is done by setting the flag to
\texttt{false} at Line~\ref{Alg1:FlagFalse}. Finally, when
$\mathcal{Q}=\emptyset$, no satisying $\mathbf{p}$ has been
found. Whether $\left[\mathbf{p}\right]$ may however contain some
satisfying $\mathbf{p}$ depends on the value of flag.

\begin{algorithm}[ht]
  \small
  \begin{algorithmic}[1]
    \protect\caption{FPS\label{alg:fps}}

    \Procedure{FPS}{$\xi_{0}$, $\xi_{u}$, $\xi_{b}$,
      $\mathbf{\left[p\right]}$,$\mathbf{\left[x\right]}$,
      $\mathbf{\left[d\right]}$} \Comment{$\xi_{0}$, $\xi_{u}$,
      $\xi_{b}$ from \eqref{eq:TestInit1}, \eqref{eq:TestUnsafe1} and
      \eqref{eq:TestBorder1}}

    \State queue $\mathcal{Q}$ := $\mathbf{\left[p\right]}$

    \State flag := \texttt{true}

    \While{$\mathcal{Q}\neq{\emptyset}$ } \State
    $\mathbf{\left[p\right]}_{0}$:= dequeue($\mathcal{Q}$)
    \Comment{Reduction of the parameter space using
      \eqref{eq:conj_contracteur2}, \eqref{eq:Contraction-tau},
      \eqref{eq:contract-params}}

    \State $\mathbf{\left[p\right]}_{0}':=\mathcal{C}_{\xi_{0}\land
      \xi_{u}}(\textrm{m}(\left[\mathbf{x}\right]),\mathbf{\left[p\right]}_{0})$
    \Comment{When $\mathbf{\left[p\right]}_{0}'=\emptyset$, there is
      no satisfying $\mathbf{p}\in\mathbf{\left[p\right]}_{0}$
    }\label{Alg1:Reduction}

    \If{$\mathbf{\left[p\right]}_{0}'=\emptyset$}
    \State continue
    \EndIf

    \State $r_{0}$:=
    CSC($\xi_{0}$,$\mathbf{\left[p\right]}_{0}',\mathbf{\left[x\right]},$$\mathbf{\left[d\right]}$)
    \Comment{Call CSC for each constraint}

    \State $r_{u}$:=
    CSC($\xi_{u}$,$\mathbf{\left[p\right]}_{0}',\mathbf{\left[x\right]},$$\mathbf{\left[d\right]}$)

    \State $r_{b}$:=
    CSC($\xi_{b}$,$\mathbf{\left[p\right]}_{0}',\mathbf{\left[x\right]},$$\mathbf{\left[d\right]}$)

    \If{($r_{0}$=\texttt{true})$\land$($r_{u}$=\texttt{true})$\land$($r_{b}$=\texttt{true})}
    \State
    return(\texttt{true},$\textrm{m}(\mathbf{\left[p\right]}_{0}')$\Comment{$\text{m}(\mathbf{\left[p\right]}_{0}')$
      is satisfying if CSCs hold true}
    \label{alg1:Alltrue}
    \EndIf

    \If{($r_{0}$\texttt{$=$false})$\vee$($r_{u}$\texttt{$=$false})$\vee$($r_{b}$\texttt{$=$false})}
    \State continue \label{Alg1:OneFalse}\Comment{no solution in
      $[\mathbf{p}]_{0}'$ if one constraint does not hold true }
    \EndIf

    \If{$\text{w}(\mathbf{\left[p\right]}_{0}')\leq\varepsilon_{\text{p}}$}
    \State flag:=\texttt{false} \Comment{no decision could be made for
      $[\mathbf{p}]_{0}'$ and thus bisect} \Else \State
    $(\mathbf{\left[p\right]}_{1},\mathbf{\left[p\right]}_{2})$:=bisect$(\mathbf{\left[p\right]}')$
    \label{Alg1:FlagFalse}

    \State enqueue($\mathbf{\left[p\right]}_{1}$) in $\mathcal{Q}$

    \State enqueue($\mathbf{\left[p\right]}_{2}$) in $\mathcal{Q}$

    \EndIf
    \EndWhile

    \If{flag=\texttt{false}} \State
    return(\texttt{unknown},$\emptyset$)\Comment{precision reached no
      conclusion can be made for $\left[\mathbf{p}\right]$}

    \EndIf

    \State return(\texttt{false},$\emptyset$) \Comment{ no valid solution in $\left[\mathbf{p}\right]$ }
    \EndProcedure
  \end{algorithmic}
\end{algorithm}

CSC, described in Algorithm~\ref{alg:csc}, determines either whether
$\tau\left(\mathbf{x},\text{m}([\mathbf{p}]),\mathbf{d}\right)$
introduced in \eqref{eq:ElementaryTest} holds true for all
$\mathbf{x}\in\left[\mathbf{x}\right]$ and all
$\mathbf{d}\in\left[\mathbf{d}\right]$ or whether there is no
$\mathbf{p}\in\left[\mathbf{p}\right]$ satisfying
\eqref{eq:ElementaryTest} for all
$\mathbf{x}\in\left[\mathbf{x}\right]$ and all
$\mathbf{d}\in\left[\mathbf{d}\right]$.  For that purpose, due to the
pessimism of inclusion functions, it may be necessary to bisect
$\left[\mathbf{x}\right]\times\left[\mathbf{d}\right]$ in subboxes
stored in a stack $\mathcal{S}$ initialized with
$\left[\mathbf{x}\right]\times\left[\mathbf{d}\right]$.  For each
subbox
$\left[\mathbf{x}\right]_{0}\times\left[\mathbf{d}\right]_{0}\subseteq\left[\mathbf{x}\right]\times\left[\mathbf{d}\right]$,
CSC determines at Line~\ref{Alg2:Validation} whether
$\textrm{m}(\mathbf{\left[p\right]}_{0})$ is satisfying. CSC tries to
prove that $\mathbf{\left[p\right]}_{0}$ does not contain any
satisfying $\mathbf{p}$. This is done at Line~\ref{Alg2:Invalidation},
for some
$\left(\mathbf{x},\mathbf{d}\right)\in\mathbf{\left[x\right]}_{0}\times\left[\mathbf{d}\right]_{0}$,
here taken as the midpoint of
$\mathbf{\left[x\right]}_{0}\times\left[\mathbf{d}\right]_{0}$, by an
evaluation of $\tau$. This is done at Line~\ref{Alg2:ContractionTau}
using the result of a contractor, as described in
\eqref{eq:contract-states} and \eqref{eq:contract-disturb}. When one
is not able to conclude and provided that
$\text{w}(\mathbf{\left[x\right]}_{0}\times\left[\mathbf{d}\right]_{0})$
is larger than $\varepsilon_{\text{x}}$, some parts of
$\mathbf{\left[x\right]}_{0}\times\left[\mathbf{d}\right]_{0}$ for
which $\textrm{m}(\mathbf{\left[p\right]}_{0})$ is satisfying are
removed at Line~\ref{Alg2:ContracteurTauBar}, before performing a
bisection and storing the resulting subboxes in $\mathcal{S}$.  When
$\text{w}(\mathbf{\left[x\right]}_{0}\times\left[\mathbf{d}\right]_{0})$
is less than $\varepsilon_{\text{x}}$ it is no more considered. As a
consequence, one is not able to determine whether
$\textrm{m}(\mathbf{\left[p\right]}_{0})$ is satisfying for all
$\left(\mathbf{x},\mathbf{d}\right)\in\mathbf{\left[x\right]}_{0}\times\left[\mathbf{d}\right]_{0}$.
One may still prove that one $\mathbf{\left[p\right]}_{0}$ does not
contain any satisfying $\mathbf{p}$ considering other subboxes of
$\left[\mathbf{x}\right]\times\left[\mathbf{d}\right]$, but not that
$\textrm{m}(\mathbf{\left[p\right]}_{0})$ is satisfying for all
$\left(\mathbf{x},\mathbf{d}\right)\in\mathbf{\left[x\right]}\times\left[\mathbf{d}\right]$.
This is indicated by setting flag to \texttt{unknown} at
Line~\ref{Alg2:Unknown}.

\begin{algorithm}[ht]
\protect\caption{CSC\label{alg:csc} }
\small
\begin{algorithmic}[1]
\Procedure{CSC}{$\tau$,$\mathbf{\left[p\right]}_{0}$,$\mathbf{\left[x\right]}$,$\left[\mathbf{d}\right]$}
\Comment{ $\tau$ is of the form \eqref{eq:ElementaryTest}}


 \State stack $S:=\mathbf{\mathbf{\left[x\right]}}\times\mathbf{\left[d\right]}$

\State flag:=\texttt{true}

\While{$S\neq\emptyset$}

\State $\mathbf{\left[x\right]}_{0}\times\left[\mathbf{d}\right]_{0}$ :=pop($S$)

\If{$[u](\mathbf{\left[x\right]}_{0},\,\textrm{m}(\mathbf{\left[p\right]}_{0}))\subseteq\mathcal{A}\lor$$[v](\mathbf{\left[x\right]}_{0},\,\textrm{m}(\mathbf{\left[p\right]}_{0}),\,\mathbf{\left[d\right]}_{0})\subseteq\mathcal{B}$}
\State continue \label{Alg2:Validation} \Comment{Validation using \eqref{eq:Validation}}
\EndIf

\If{$[u]\left(\text{m}\left(\mathbf{\left[x\right]}_{0}\right),\mathbf{\left[p\right]}_{0}\right)\cap\mathcal{A}=\emptyset\land[v]\left(\text{m}\left(\mathbf{\left[x\right]}_{0}\right),\mathbf{\left[p\right]}_{0},\text{m}\left(\mathbf{\left[d\right]}_{0}\right)\right)\cap\mathcal{B}=\emptyset$}
\State return(\texttt{false}) \Comment{Reduction of the parameter
  space using \eqref{eq:invalidation}}
 \label{Alg2:Invalidation}
\EndIf

\State$(\mathbf{\left[x\right]}_{0}',\mathbf{\left[p\right]}_{0}',\mathbf{\left[d\right]}_{0}'):=\mathcal{C}_{\tau}(\mathbf{\left[x\right]}_{0},\mathbf{\left[p\right]}_{0},\mathbf{\left[d\right]}_{0})$\Comment{Reduction
  of the parameter space applying \eqref{eq:Contraction-tau}}
\label{Alg2:ContractionTau}

\If{$\mathbf{\left[x\right]}_{0}'\neq\mathbf{\left[x\right]}_{0}\lor\mathbf{\left[d\right]}_{0}'\neq\mathbf{\left[d\right]}_{0}$}
\State return(\texttt{false}) \Comment{Consequence of
  \eqref{eq:contract-states}, \eqref{eq:contract-disturb}}
\EndIf

\If{($\text{w}(\mathbf{\left[x\right]}_{0}\times\left[\mathbf{d}\right]_{0})\leq\varepsilon_{\text{x}}$)}
\State flag:=\texttt{unknown} \label{Alg2:Unknown}
\Else
\State $(\mathbf{\left[x\right]}_{0}'',\mathbf{\left[d\right]}_{0}''):=\mathcal{C}_{\overline{\tau}}(\mathbf{\left[x\right]}_{0},\textrm{m}(\mathbf{\left[p\right]}_{0}),\mathbf{\left[d\right]}_{0})$ \label{Alg2:ContracteurTauBar}\Comment{Reduction of the state space using  \eqref{eq:contraction-state-space}, \eqref{eq:contract-neg} and bisection}

\State $(\mathbf{\left[x\right]}_{1}\ensuremath{\times\mathbf{\left[d\right]}_{1}},\mathbf{\left[x\right]}_{2}\times\mathbf{\left[d\right]}_{2})$:=bisect($\mathbf{\left[x\right]}_{0}''\times\mathbf{\left[d\right]}_{0}''$)

\State push($\mathbf{\left[x\right]}_{1}\ensuremath{\times\mathbf{\left[d\right]}_{1}}$) in $S$

\State push($\mathbf{\left[x\right]}_{2}\times\mathbf{\left[d\right]}_{2}$)in $S$
\EndIf
\EndWhile

\State return(flag)
\EndProcedure
\end{algorithmic}
\end{algorithm}

\section{Examples}

\label{sec:example} This section presents experiments for the characterization of barrier
functions. The considered dynamical systems are described first before
providing numerical results, comparison of different approaches and
discussions.

\subsection{Dynamical system descriptions}

For the following examples, one provides the dynamics of the system,
the constraints $g_{0}$ and $g_{u}$ for the definition of the sets
$\mathcal{X}_{0}$ and $\mathcal{X}_{u}$, the state space
$[\mathbf{x}]$, and the parametric expression of the barrier
function. In all cases, the parameter space is chosen as
$\mathbf{\left[p\right]}=[-10,10]^{m}$ where $m$ is the number of
parameters.

\begin{example}
  \label{exple:p0}
  Consider the system
  \[
  \begin{pmatrix}\dot{x}_{1}\\
    \dot{x}_{2}
  \end{pmatrix}=\begin{pmatrix}x_{1}+x_{2}\\
    x_{1}x_{2}-0.5x_{2}^{2}
  \end{pmatrix}
  \]
  with $g_{0}(\mathbf{x})=(x_{1}+1.25)^{2}+(x_{2}-1.25)^{2}-0.05$,
  $g_{u}(\mathbf{x})=(x_{1}+2.5)^{2}+(x_{2}-0.8)^{2}-0.05$, and $\mathbf{\left[x\right]}=[-10^{3},0]\times[-10^{3},10^{3}]$.
  The parametric barrier function is $B(\mathbf{x},\mathbf{p})=\frac{p_{1}p_{2}(x_{0}+p_{3})}{(x_{0}+p_{3})^{2}+p_{2}^{2}}+x_{1}+p_{4}$.\end{example}

\begin{example}
\label{exple:parillo}
  Consider the system from~\cite{ahmadi2011globally}
  \[
  \begin{pmatrix}\dot{x}_{1}\\
    \dot{x}_{2}
  \end{pmatrix}=\begin{pmatrix}-x_{1}+x_{1}x_{2}\\
    -x_{2}
  \end{pmatrix}
  \]
  with $g_{0}(\mathbf{x})=(x_{1}-1.125)^{2}+(x_{2}-0.625)^{2}-0.0125$,
  $g_{u}(\mathbf{x})=(x_{1}-0.875)^{2}+(x_{2}-0.125)^{2}-0.0125$, and
  $\left[x\right]=[-100,100]\times[-100,100]$. The parametric barrier
  function used is
  $B(\mathbf{x},\mathbf{p})=\ln(p_{1}x_{1})-\ln(x_{2})+p_{2}x_{2}+p_{3}$.
\end{example}

\begin{example}
  \label{exple:saturation}
  Consider the system from~\cite{papachristodoulou2005analysis}

  \[
  \begin{pmatrix}\dot{x}_{1}\\
    \dot{x}_{2}
  \end{pmatrix}=\begin{pmatrix}x_{2}\\
    -\frac{x1+x2}{\sqrt{1+(x1+x2)^{2}}}
  \end{pmatrix}
  \]
  with $g_{0}(\mathbf{x})=x_{1}{}^{2}+x_{2}{}^{2}-0.5$, $g_{u}(\mathbf{x})=(x_{1}-3.5)^{2}+(x_{2}-0.5)^{2}-0.5$
  , and$\mathbf{\left[x\right]}=[-10^{3},10^{3}]\times[-100,100]$.
  A quadratic parametric barrier function is chosen $B(\mathbf{x},\mathbf{p})=p_{1}x_{1}^{2}+p_{2}x_{2}^{2}+p_{3}x_{1}x_{2}+p_{4}x_{1}+p_{5}x_{2}+p_{6}$.
\end{example}

\begin{example}
  \label{exple:p3}
  Consider the disturbed system from
  \cite{prajna2004safety}
  \[
  \begin{pmatrix}\dot{x}_{1}\\
    \dot{x}_{2}
  \end{pmatrix}=\begin{pmatrix}x_{2}\\
    -x_{1}+\frac{d}{3}x_{1}^{3}-x_{2}
  \end{pmatrix}
  \]
  with $g_{0}(\mathbf{x})=(x_{1}-1.5)^{2}+x_{2}{}^{2}-0.25$, $g_{u}(\mathbf{x})=(x_{1}+0.8)^{2}+(x_{2}+1)^{2}-0.25$,
  $\mathbf{\left[x\right]}=[-100,100]\times[-10,10]$, and $d\in[0.9,1.1]$.
  The parametric barrier function $B(\mathbf{x},\mathbf{p})=p_{1}x_{1}^{2}+p_{2}x_{2}^{2}+p_{3}x_{1}x_{2}+p_{4}x_{1}+p_{5}x_{2}+p_{6}$
  is considered.
\end{example}

\begin{example}
  \label{exple:p9}
  Consider the system with a limit cycle
  \[
  \begin{pmatrix}\dot{x}_{1}\\
    \dot{x}_{2}
  \end{pmatrix}=\begin{pmatrix}x_{2}+(1-x_{1}^{2}-x_{2}^{2})x_{1}+\ln(x_{1}^{2}+1)\\
    -x_{1}+(1-x_{1}^{2}-x_{2}^{2})x_{2}+\ln(x_{2}^{2}+1)
  \end{pmatrix}
  \]
  with $g_{0}(\mathbf{x})=(x_{1}-1)^{2}+(x_{2}+1.5)^{2}-0.05$, $g_{u}(\mathbf{x})=(x_{1}+0.6)+(x_{2}-1)^{2}-0.05$,
  and $\mathbf{\left[x\right]}=[-10^{3},10^{3}]\times[-10^{3},10^{3}]$.
  The parametric barrier function used is $B(\mathbf{x},\mathbf{p})=\left(\frac{x_{1}+p_{1}}{p_{2}}\right)^{2}+\left(\frac{x_{2}+p_{3}}{p_{4}}\right)^{2}-1$.

\end{example}

\begin{example}
  \label{exple:lorenz}
  Consider the Lorenz system from \cite{vanvevcek1996control}
  \[
  \begin{pmatrix}\dot{x_{1}}\\
    \dot{x_{2}}\\
    \dot{x_{3}}
  \end{pmatrix}=\begin{pmatrix}10(x_{2}-x_{1})\\
    x_{1}(28-x_{3})-x_{2}\\
    x_{1}x_{2}-\frac{8}{3}x_{3}
  \end{pmatrix}
  \]
  with $g_{0}(\mathbf{x})=(x_{1}+14.5)^{2}+(x_{2}+14.5){}^{2}+(x_{3}-12.5)^{2}-0.25$,
  $g_{u}(\mathbf{x})=(x_{1}+16.5)^{2}+(x_{2}+14.5){}^{2}+(x_{3}-2.5)^{2}-0.25$,
  and $\mathbf{\left[x\right]}=[-20,20]\times[-20,0]\times[-20,20]$.
  The considered parametric barrier function is $B(\mathbf{x},\mathbf{p})=p_{1}x_{1}^{2}+p_{2}x_{1}+p_{3}x_{3}+p_{4}$.

\end{example}

\begin{example}
  \label{exple:sixdim}
  Consider the system from \cite{papachristodoulou2002construction}
  \[
  \begin{pmatrix}\dot{x_{1}}\\
    \dot{x_{2}}\\
    \dot{x_{3}}\\
    \dot{x_{4}}\\
    \dot{x_{5}}\\
    \dot{x_{6}}
  \end{pmatrix}=\begin{pmatrix}-x_{1}+4x_{2}-6x_{3}x_{4}\\
    -x_{1}-x_{2}+x_{5}^{3}\\
    x_{1}x_{4}-x_{3}+x_{4}x_{6}\\
    x_{1}x_{3}+x_{3}x_{6}-x_{4}^{3}\\
    -2x_{2}^{3}-x_{5}+x_{6}\\
    -3x_{3}x_{4}-x_{5}^{3}-x_{6}
  \end{pmatrix}
  \]
  with
  $g_{0}(\mathbf{x})=(x_{1}-3.05)^{2}+(x_{2}-3.05){}^{2}+(x_{3}-3.05)^{2}+(x_{4}-3.05)^{2}+(x_{5}-3.05)^{2}+(x_{6}-3.05)^{2}-0.0001$,
  $g_{u}(\mathbf{x})=(x_{1}-7.05)^{2}+(x_{2}-3.05){}^{2}+(x_{3}-7.05)^{2}+(x_{4}-7.05)^{2}+(x_{5}-7.05)^{2}+(x_{6}-7.05)^{2}-0.0001$,
  and
  $\mathbf{\left[x\right]}=[0,10]\times[0,10]\times[2,10]\times[0,10]\times[0,10]\times[0,10]$.
  The considered parametric barrier function is
  $B(\mathbf{x},\mathbf{p})=p_{1}x_{1}^{2}+p_{2}x_{2}^{4}+p_{3}x_{3}^{2}+p_{4}x_{4}^{2}+p_{5}x_{5}^{4}+p_{6}x_{6}^{2}+p_{7}$.

\end{example}

\subsection{Experimental conditions and results}

CSC-FPS, presented in Section~\ref{sec:algorithm}, has been
implemented using the IBEX
library~\cite{araya2012upper,chabert2009contractor}.  The selection of
candidate barrier functions is performed choosing polynomials with
increasing degree, except for Examples~\ref{exple:p0}, \ref{exple:parillo},
and \ref{exple:p3}, where parametric
functions taken from \cite{mathcurve,ahmadi2011globally} are
considered.

\begin{figure}
  \centering
  \includegraphics{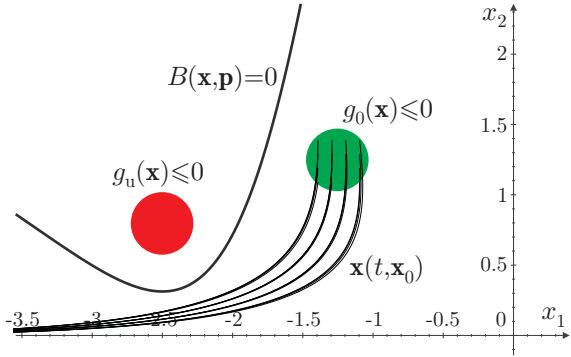}

  \protect\caption{Results for Example~\ref{exple:p0}.
    \label{P0-fig}}
\end{figure}

\begin{figure}
  \centering
  \includegraphics{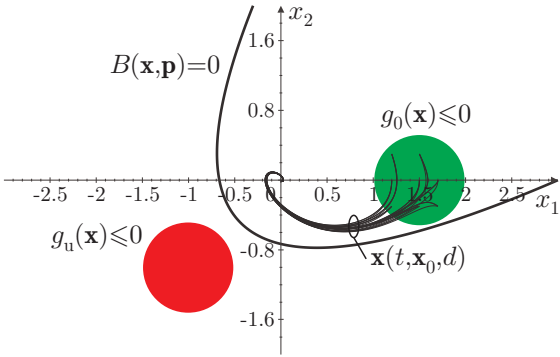}

  \protect\caption{Results for Example~\ref{exple:p3}
    with various values of the disturbance
    \textbf{$d$}.
    \label{P3-fig}}
\end{figure}

\begin{figure}
  \centering
  \includegraphics{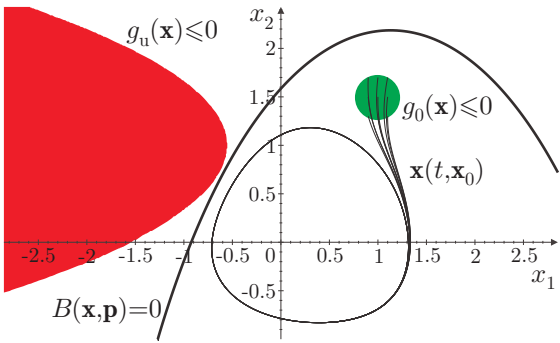}
  \protect\caption{Results for Example~\ref{exple:p9}. \label{fig:p9}}
\end{figure}

For each example, the computing time to get a valid barrier function
and the number of bisections of the search box
$\left[\mathbf{p}\right]$ are
provided. Table~\ref{tab:no-and-contract} summarizes the results for
the versions of CSC-FPS with and without contractors. As in
\cite{Jaulin1996guaranteed}, we choose
$\varepsilon_{\text{x}}=10^{-1}$ and
$\varepsilon_{\text{p}}=10^{-5}$. Computations were done on an Intel
core $1.7$~GHz processor with $8$~GB of RAM. If after $30$ minutes of
computations no valid barrier function has been found, the search is
stopped. This is denoted by T.O. (time out) in Table~\ref{tab:no-and-contract}. Moreover,
for Examples~\ref{exple:p0}, \ref{exple:p3}, and \ref{exple:p9},
graphical representation of the computed barrier functions are
provided. In Figures~\ref{P0-fig}, \ref{P3-fig}, and \ref{fig:p9},
$\mathcal{X}_{0}$ is in green, $\mathcal{X}_{u}$ is in red, the bold
line is the barrier function and some trajectories starting from
$\mathcal{X}_{0}$ are also represented.

The results in Table~\ref{tab:no-and-contract} show the importance of
contractors which are beneficial in all cases.  Thanks to contractors,
valid barrier functions were obtained for all examples, which is not
the case employing the original version of CSC-FPS proposed in
\cite{Jaulin1996guaranteed}. In theory, both FPS and CSC are of
exponential complexity in the dimension $m$ of the parameter space and
$n$ of the state space. In practice, contractors allow, on the
considered examples, to break this complexity and to consider
high-dimensional problems.

The search for barrier functions using the relaxed version
\eqref{eq:Tborderrelaxed} of the constraint \eqref{eq:TestBorder} as
in \cite{prajna2004safety}, has also been performed using the version
of our algorithms with contractors considering the same parametric
barrier functions. We were not able to find a valid barrier function
in less than $30$ minutes. This shows the detrimental effect of the
relaxation \eqref{eq:Tborderrelaxed} on the search technique.

Some examples were also addressed using RSolver, which is a tool to
solve some classes of QCSP \cite{ratschan2006efficient}. Nevertheless,
this requires some modifications of the examples, since RSolver does
not address dynamics or constraints involving divisions
\cite{Rsolverman}.  Moreover, the type of constraints processed by
RSolver, for problems such as \eqref{eq:tau-QSCP}, allows only
parametric barrier functions that are linear in the
parameters. Outer-approximating boxes for the initial and the unsafe
regions $\mathcal{X}_{0}$, $\mathcal{X}_{u}$ defined by $g_{0}$ and
$g_{u}$ are given to Rsolver. As a consequence,
Examples~\ref{exple:p0}, \ref{exple:parillo} and
\ref{exple:saturation} could not be considered. Only
Examples~\ref{exple:p3}, \ref{exple:p9}, and \ref{exple:lorenz} were
tested by Rsolver, which was able to find a satisfying barrier
function only for Lorenz in less than $1$s. RSolver was unable to find
a solution for Examples~\ref{exple:p3} and \ref{exple:p9}. Rsolver is
well designed for problems with parametric barriers linear in the
parameters, but has difficulties for non-linear problem and non-convex
constraints such as~\eqref{eq:TestBorder}.




\begin{table}[t]
  \centering
  \caption{Results of CSC-FPS without and with contractors}
  \begin{tabular}{|c|c|c|c|c|c|c|}
    \hline
    \multicolumn{3}{| c}{} & \multicolumn{2}{| c|}{Without contr.} &
    \multicolumn{2}{c|}{With contr.}
    \tabularnewline
    \hline
    Example  & $n$ & $m$ & time  &  bisect. & time & bisect.  \tabularnewline
    \hline
    1  & 2 & 4 & $36$s & $4520$ & $16$s & $4553$ \tabularnewline
    \hline
    2  & 2 & 3 & T.O. & / & $1$s & $159$\tabularnewline
    \hline
    3  & 2 & 6 & $1133$s & $20388$ & $1$s & $6$ \tabularnewline
    \hline
    4  & 2 & 6 & $253$s  & $14733$ & $7$s & $435$ \tabularnewline
    \hline
    5  & 2 & 4 & T.O. & $/$ & $98$s & $4072$ \tabularnewline
    \hline
    6  & 3 & 4 & $167$s  & $1753$ & $21$s & $47$ \tabularnewline
    \hline
    7 & 6 & 7 & $697$s & $67600$ & $1$s & $261$ \tabularnewline
    \hline
  \end{tabular}
  \label{tab:no-and-contract}
\end{table}

\section{Conclusion}

\label{sec:conclusion} This paper presents a new method to find parametric barrier functions
for nonlinear continuous-time perturbed dynamical systems. The proposed
technique has no restriction regarding the dynamics nor the form of
the barrier function. The search for barrier functions is formulated
as an interval quantified constraint statisfaction problem. A branch-and-prune
algorithm proposed in \cite{Jaulin1996guaranteed} has been supplemented
with contractors to address this problem. Contractors are instrumental
in solving problems with large number of parameters. The proposed
approach can thus find barrier functions for a large class of possibly
perturbed dynamical systems.

Alternative techniques based on RSolver may be significantly more
efficient for some specific classes of problems where the parameters
appear linearly in the parametric barrier functions. A combination
of RSolver and our approach may be useful to improve the global efficiency
of barrier function caracterization. 

Future work will be dedicated to the search for the class of parametric
barrier functions to consider. This may be done by exploring a library
of candidate barrier functions. In our approach rejection of a candidate
function occurs mainly after a timeout. Even if contractors aiming
at eliminating some parts of the parameter space were defined, their
efficiency is limited. Better contractors for that purpose may be
very helpful.

An other future research direction is to extend the proposed method
to hybrid systems as done in~\cite{prajna2004safety}, \emph{i.e.},
to consider a set of quantified constraints for each location of an
hybrid automaton and the constraints associated to the transitions.

\bibliographystyle{plain}
\bibliography{paper}

\end{document}